\documentclass{article}

\usepackage{arxiv}

\usepackage[utf8]{inputenc} % allow utf-8 input
\usepackage[T1]{fontenc}    % use 8-bit T1 fonts
\usepackage{hyperref}       % hyperlinks
\usepackage{url}            % simple URL typesetting
\usepackage{booktabs}       % professional-quality tables
\usepackage{amsfonts}       % blackboard math symbols
\usepackage{nicefrac}       % compact symbols for 1/2, etc.
\usepackage{microtype}      % microtypography
\usepackage{lipsum}
\usepackage{graphicx}
\graphicspath{ {./Images/} }
\usepackage{tabu}
\usepackage{amsmath}
\usepackage{multirow}
\usepackage{threeparttable}
\usepackage{algorithm}
\usepackage{algorithmic}
\usepackage[title]{appendix}
\usepackage[english]{babel}
\usepackage{amsthm}

\newtheorem{theorem}{Theorem}
\newtheorem{remark}{Remark}
\newtheorem{assumption}{Assumption}

\title{On the Exactness of an Energy-efficient Train Control model based on Convex Optimization}

\author{ \href{https://orcid.org/0000-0001-5361-2463}{\includegraphics[scale=0.06]{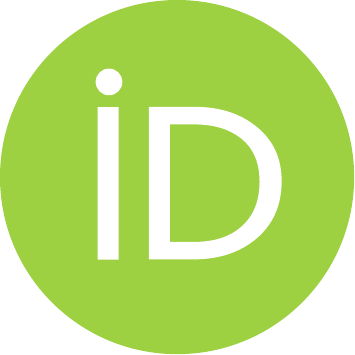}\hspace{1mm}Shaofeng~ Lu}\thanks{Corresponding author, https://lushaofeng.github.io}\\
Shien-ming Wu School of Intelligent Engineering\\
  South China University of Technology\\
  Guangzhou, China 511442 \\
	\texttt{lushaofeng@scut.edu.cn} \\
	\And
	 Minling Feng \\
  Shien-ming Wu School of Intelligent Engineering\\
  South China University of Technology\\
  Guangzhou, China 511442 \\
  \texttt{Minling.Feng@outlook.com} \\
   \And
 Kunpeng Wu \\
  Shien-ming Wu School of Intelligent Engineering\\
  South China University of Technology\\
  Guangzhou, China 511442 \\
  \texttt{Kunpeng.Wu@outlook.com} \\
}

\renewcommand{\shorttitle}{Exactness Proof of EETC Model based on Convex Optimization}

\hypersetup{
pdftitle={On the Exactness of an Energy-efficient Train Control model based on Convex Optimization},
pdfsubject={Optimization and Control},
pdfauthor={Shaofeng Lu, Minling Feng and Kunpeng Wu},
pdfkeywords={Convex Optimization, Energy-efficient Train Control, Convexification and Relaxations},
}

\begin{document}
\maketitle

\begin{abstract}
In this paper, we demonstrate the exactness proof for the energy-efficient train control (EETC) model based on convex optimization. The  proof of exactness shows that the convex optimization model will share the same optimization results with the initial model on which the convex relaxations are conducted. We first show how the relaxation on the initial non-convex model is conducted and provide analysis to show that the relaxations are convex constraints and the relaxed model is thus a convex model. Subsequently, we prove that the relaxed convex model will always achieve its optimal solution on the initial equality constraints and the optimal solution achieved by convex optimization will be the same as the one obtained by the initial non-convex model and the relaxations applied are exact. A numerical verification has been conducted based on a typical urban rail system with a steep gradient. The results of this paper shed lights on further applications of convex optimization on energy-efficient train control and relevant areas related to operation and control of low-carbon transportation systems. 
\end{abstract}

%\setcounter{tocdepth}{2}
%\tableofcontents

\section{Introduction}
Energy-efficient train control (EETC) tries to locate the most energy-efficient train speed trajectory resulted from the optimal control strategies of the train dynamic system. The past studies in the field generally adopted two main types of methods: indirect method based on the optimal theory and direct methods based on different optimization techniques.

On the one hand, indirect method obtained its name from the way of locating the optimal solution by using the Pontryagin's Maximum Principle (PMP). According to PMP, the optimal solution's necessary condition can be obtained, upon which the solution can be evaluated using the necessary condition. The optimal solution can be obtained by indirectly solving the co-state variable's differential equations with the help of the necessary conditions. On the other hand, direct methods try to set up mathematical models and apply computer algorithms to locate the optimal solution. Our proposed EETC model with its solution approach falls into the direct method and in general, we are proposing a mathematical model solved by efficient algorithms to quickly locate the optimal solutions. These solutions are generally represented by the train speed and its corresponding train control strategies, i.e. acceleration, braking and cruising etc. We refer readers to numerous influential academic papers on both indirect methods \cite{Khmelnitsky2000,Howlett2000,Liu2003,Albrecht2016a,Albrecht2016b} and direct methods\cite{Goverde2020,WangY2013,WangP2016,YeH2017,WuC2019a}. 

Based on the previous research outcomes of \cite{FengM2022convex}, in this paper, we will presents two direct EETC optimization models.  The first model, referred to as ``\textbf{Model A}'' is an non-convex model and the second mode is the relaxed ``\textbf{Model B}''. We focus on the topic how the two optimization models can share same optimization results after the relaxation has been conducted from one to another.  We first present two models, offer some analysis on the convexity of Model B in Section \ref{sec:models}, then provide the exactness proof of relaxations in Section \ref{sec:exact_proof} and demonstrate a numerical verification using a typical urban rail case with a steep downhill gradient in Section \ref{sec:num_veri}. A conclusion is drawn in Section \ref{sec:conclusion}.  

\section{Two EETC Models}
\label{sec:models}
%In Fig. \ref{fig:two_models}, we present the details of the two models. The non-convex model is shown on the left column with 7 sets of constraints (a)-(f) imposed on each distance segment. By relaxation, constraint (d) and (f) becomes convex.

\subsection{The non-convex EETC model}
A similar introduction of this model can be found in \cite{FengM2022convex}, for the sake of coherence, we still include the key modeling details in this section.  The journey is divided into $N$ segments. For the sake of simplicity, the length of each segment is equal and assumed to be $\Delta d$. The total length of journey is denoted by $D$.
We have
\begin{equation}
\label{for:jou_length_constr}
    D=\sum_{i=1}^N \Delta d_i
\end{equation}
where $i$ is the index number corresponding to each candidate speed point along the journey.

The whole journey includes $N + 1$ speed variable $v_i$ corresponding to $N$ segments. $v_0$ and $v_N$ are initial and final candidate speed. Each speed $v_i$ is
constrained by the speed limit $V_{i,lim}$ as presented by:
\begin{equation}
\label{for:speed_lim_constr}
    0<v_i\leq V_{i,lim}
\end{equation}
where $V_{i,lim}$ is a parameter setting the speed limit for each $v_i$. In the proof demonstrated later, we make a special case for $v_0$ and assume $v_0=0$ to ensure the train departs from zero speed. 

The sum of elapsed time in each segment equals to the total journey time $T$ as shown in \eqref{for:total_time_constr}. 
\begin{equation}
\label{for:total_time_constr}
    \sum_{i=1}^N\frac{\Delta d}{v_i}=T. 
\end{equation}

The drag force can be calculated by the Davis equation in \eqref{for:Davis_constr}:
\begin{equation}
\label{for:Davis_constr}
    f_{i,d}=A+Bv_i+Cv_i^2
\end{equation}
where the parameters $A$, $B$ and $C$ are the Davis coefficients. 

The principle of energy conservation during train running can be presented by \eqref{for:ene_conserv_1}:
\begin{align}
     \label{for:ene_conserv_1}
     F_i\Delta d &=0.5M(v_i^2 -v_{i-1}^2)+f_{i,d}\Delta d +Mg\Delta H_i\nonumber\\
     &=0.5M(v_i^2 -v_{i-1}^2)+(A+B v_i+Cv_i^2)\Delta d +Mg\Delta H_i
\end{align}

where $F_i$ represents the tractive or braking effort imposed the train based on the principle of conservation of energy, $M$ is the total mass of the train which can be varied to consider the rotary mass, and $g$ is acceleration rate due to gravity. $\Delta H_i$ is the altitude difference between the current distance $d_i$ and previous position $d_{i-1}$ and $\Delta H_i$ is positive when it is going uphill.

The inequality constraints are imposed to ensure that the train effort and power does not exceed its boundaries:

\begin{equation}
\label{for:force_boundary}
    -F_{max}\leq F_i\leq F_{max}
\end{equation}

\begin{equation}
\label{for:power_boundary}
    \frac{-P_{b,max}}{v_i}\leq F_i\leq \frac{P_{t,max}}{v_i}
\end{equation}

where $P_{b,max}$ is the maximum braking power and $P_{t,max}$ is the maximum traction power.  \eqref{for:force_boundary} is to denote the force boundary condition and \eqref{for:power_boundary} is for the power boundary. 

The electrical energy consumed during traction and generated during the regenerative braking can be modeled by:
\begin{subequations}
\label{for:elec_energy}
\begin{align}
    E_i&\geq F_i\Delta d/\eta_t\label{subfor:elec_energy_a}\\
    E_i&\geq F_i\Delta d \eta_b\label{subfor:elec_energy_b}
\end{align}
\end{subequations}

where $1\geq \eta_t>0$ and $1\geq \eta_b>0$ are the motor efficiency during traction and braking respectively. When $F_i$ is positive \eqref{subfor:elec_energy_b} is relaxed and otherwise \eqref{subfor:elec_energy_a} is relaxed. $\eta_b$ can be set to sufficiently small number so that no regenerative braking energy will be considered.

The total net energy is defined as the total positive tractive energy plus the total negative regenerative braking energy. 

The non-convex EETC model can be presented by:

\begin{align}
\label{for:EETC_noncvx}
    &\mathrm{minimize}~~f_{\rm obj}(E)=\sum_{i=1}^N E_i \\
    &\mathrm{Subject~to:~}\eqref{for:jou_length_constr}-\eqref{for:elec_energy}\nonumber
\end{align}

Since the objective function is to minimize, $E_i$ will always approach $F_i \Delta d \eta_t$ or $F_i\Delta d \eta_b$ as shown in \eqref{for:elec_energy}.

\subsection{The relaxed EETC model}
In the relaxed convex model, two new sets of variables are introduced: $\alpha_i$ and $\beta_i$. 

%We refer these two sets of variables as the ``auxiliary variables'' and $1/v_i$ and $v_i^2$ as the ``original variables''. 

We add two constraints defined in \eqref{for:alpha_i} and \eqref{for:beta_i}. 

\begin{equation}
\label{for:alpha_i}
    \alpha_i\geq 1/v_i \mathrm{~~for~~i= 1, 2, \cdots,N}
\end{equation}

\begin{equation}
\label{for:beta_i}
    \beta_i\geq v_i^2 \mathrm{~~for~~i= 1, 2, \cdots,N-1}
\end{equation}

The constraint \eqref{for:speed_lim_constr} is transformed into \eqref{for:speed_lim_constr_cvx}:
\begin{equation}
    \label{for:speed_lim_constr_cvx}
    \beta_i\leq V_{i,lim}^2.
\end{equation}

The constraint \eqref{for:total_time_constr} is transformed into \eqref{for:total_time_constr_cvx}:

\begin{equation}
\label{for:total_time_constr_cvx}
    \sum_{i=1}^N\Delta d\alpha_i=T. 
\end{equation}

The constraint \eqref{for:Davis_constr} is transformed into \eqref{for:Davis_constr_cvx}:

\begin{equation}
\label{for:Davis_constr_cvx}
    f_{i,d}=A+Bv_i+C\beta_i. 
\end{equation}

The constraint \eqref{for:ene_conserv_1} is transformed into \eqref{for:ene_conserv_1_cvx}:

\begin{align}
     \label{for:ene_conserv_1_cvx}
     F_i\Delta d &=0.5M(\beta_i -\beta_{i-1})+f_{i,d}\Delta d +Mg\Delta H_i\nonumber\\
     &=0.5M(\beta_i -\beta_{i-1})+(A+B v_i+C\beta_i)\Delta d +Mg\Delta H_i.
\end{align}

The power boundary constraint \eqref{for:power_boundary} is transformed into \eqref{for:power_boundary_cvx}:

\begin{equation}
\label{for:power_boundary_cvx}
    -P_{b,max} \alpha_i\leq F_i\leq P_{t,max} \alpha_i.
\end{equation}

The relaxed EETC model can be presented by \eqref{for:EETC_cvx}:

\begin{align}
\label{for:EETC_cvx}
    &\mathrm{minimize}~~f_{\rm obj}(E)=\sum_{i=1}^N E_i \\
    &\mathrm{Subject~to:~}\eqref{for:jou_length_constr},\eqref{for:force_boundary}, \eqref{for:elec_energy}, \eqref{for:alpha_i}- \eqref{for:power_boundary_cvx}. \nonumber
\end{align}

The relaxed EETC model is built with two relaxation constraints defined in \eqref{for:alpha_i}  and \eqref{for:beta_i}. In Subsection \ref{subsec:convexity_analysis}, we demonstrate that these two constraints are convex and the new EETC defined in \eqref{for:EETC_cvx} is a convex optimization model. In Section \ref{sec:exact_proof}, we prove that the two relaxation constraints in \eqref{for:alpha_i}  and \eqref{for:beta_i} are exact. In other words, the optimal solution will always be attained on the equality of the constraints. 

\subsection{Analysis of convexity of two relaxed constraints}
\label{subsec:convexity_analysis}
The constraint \eqref{for:alpha_i} can be easily transformed into a Second-order Conic Program (SOCP) given that both $\alpha_i \geq 0$ and $v_i\geq 0$. 

We conduct the following transformation on \eqref{for:alpha_i} as presented in \eqref{for:SOCP_1}:

\begin{subequations}
\label{for:SOCP_1}
\begin{align}
1 &\leq \alpha_i v_i \\
2^2+(\alpha_i-v_i)^2&\leq \alpha_i^2 +v_i^2+2\alpha_i v_i\\
\sqrt{2^2+(\alpha_i-v_i)^2}&\leq \alpha_i +v_i.
\end{align}
\end{subequations}

Define new variables $x_i$, $y_i$ and $z_i$ in \eqref{for:SOCP_2}:
\begin{subequations}
\label{for:SOCP_2}
\begin{align}
    x_i&=2\\
    y_i&=\alpha_i-v_i\\
    z_i&=\alpha_i+v_i.
\end{align}
\end{subequations}

This transforms the program \eqref{for:SOCP_1} into a standard SOCP as shown in \eqref{for:SOCP_3}:
\begin{align}
    \label{for:SOCP_3}
    \sqrt{x_i^2+y_i^2}\leq z_i. 
\end{align}

The definition in \eqref{for:SOCP_2} represents the intersection of a second-order cone with affine sets and this does not change the convexity of the constraint in \eqref{for:SOCP_1}. It thus verifies that constraint \eqref{for:alpha_i} is a convex SOCP constraint given that $v_i> 0$ and $\alpha_i> 0$. 

On the other hand, \eqref{for:beta_i} is an inequality involving quadratic term $v_i^2$ and it is a convex constraint which can be verified by the first-order condition of convex function \cite[Page 83]{Convex_boyd}. For the sake of simplicity, we consider the following function to represent $v_i^2\leq \beta_i$:
\begin{equation}
    \label{for:beta_i_f}
    f(\mathbf{x})=x_1^2-x_2~~~ x_1\geq 0~~\mathrm{and}~~x_2\geq 0
\end{equation}
where $\mathbf{x}=[x_{1},x_{2}]^T$ is to represent the decision vector $[v_i, \beta_i]^T$.

The first-order condition of convexity states that if function $f$ is differentiable, then $f$ is convex if and only if the feasible domain of $f$ is convex and
\begin{equation}
\label{for:1st_order_condition}
    f(y)\geq f(x)+\nabla f(x)^T(y-x).
\end{equation}

The feasible domain of function $f(\mathbf{x})$ is convex which can be proved by the convex set definition. We conduct the following deduction on \eqref{for:beta_i_f} based on the first-order conditions for convexity:
\begin{align}
    &f_i(\mathbf{y_i})-[f_i(\mathbf{x_i})+\nabla f(x)^T(\mathbf{y_i}-\mathbf{x_i})]\nonumber\\
    &=y_1^2-y_2-[x_1^2-x_2+2x_1(y_1-x_1)-(y_2-x_2)]\nonumber\\
    &=y_1^2-y_2-(-y_2-x_1^2+2x_1y_1)\nonumber\\
    &=(y_1-x_1)^2\geq 0
\end{align}
This shows that the constraint $v_i^2\leq \beta_i$ is a convex constraint. It can be seen that all constraints in Model B is either affine equality or convex inequality. The objective function \eqref{for:EETC_cvx} is also an affine function. Optimization problem based on Model B is a convex optimization problem. 

\section{Proof of Exactness of Relaxations}
\label{sec:exact_proof}

We use Model A to refer to the model defined by \eqref{for:EETC_noncvx} and Model B to refer to the model defined by \eqref{for:EETC_cvx}. 

\begin{theorem}
\label{the:exactness}
The optimal solution of Model A is also the optimal solution of Model B and thus the relaxation constraints \eqref{for:alpha_i} and \eqref{for:beta_i} in Model B are exact.
\end{theorem}

Assume the optimal solution for the model defined by Model A is $\mathcal{S}_1^*$ and the optimal solution for Model B is $\mathcal{S}_2^*$. 

Since Model B is a relaxed model of Model A, $f_{\mathrm{obj}}(\mathcal{S}_2^*)\leq f_{\mathrm{obj}}(\mathcal{S}_1^*)$. On the other hand, if the optimal solution of Model B is always a feasible solution of Model A, it can be deducted that $f_{\mathrm{obj}}(\mathcal{S}_1^*)\leq f_{\mathrm{obj}}(\mathcal{S}_2^*)$. This assumption is valid as long as the optimal solution can always be attained on the equality of \eqref{for:alpha_i} and \eqref{for:beta_i}. As a result, when the optimal solution is achieved on the equality of \eqref{for:alpha_i} and \eqref{for:beta_i}, we have $f_{\mathrm{obj}}(\mathcal{S}_1^*)=f_{\mathrm{obj}}(\mathcal{S}_2^*)$. 

Therefore, in order to prove Theorem \ref{the:exactness}, it suffices to show that any optimal solution of the model defined by \eqref{for:EETC_cvx} attains equality in \eqref{for:alpha_i} and \eqref{for:beta_i}.

We will prove that equality will be attained in \eqref{for:alpha_i} for any optimal solution in Subsection \ref{subsec:PoE_part1} and prove the equality is also attained in \eqref{for:beta_i} for any optimal solutions in Subsection \ref{subsec:PoE_part2}. 

We make the following assumption:
\begin{assumption}
\label{assum:1}
We assume that $Mg\Delta H_1=0$ and $\beta_0=0$. It derives:
\begin{equation}
    F_1\Delta d=0.5M\beta_1+\Delta d(A+Bv_1+c\beta_1). 
\end{equation}
\end{assumption}

Based on Assumption \ref{assum:1}, it is known that $F_1\Delta d$ is a positive number and its reduction caused by a positive value $\beta_1$ and $v_1$ will not violate the force and power boundary as defined in \eqref{for:force_boundary} and \eqref{for:power_boundary_cvx}. It is known that $v_0=\beta_0=0$ and it is also a reasonable assumption that the train will start from a zero speed. 

\begin{assumption}
\label{assum:2}
We assume that a strict inequality exists for $v_N^2<\beta_N$. 
\end{assumption}

Assumption \ref{assum:2} is used to demonstrate that the model will always find a feasible solution to keep \eqref{for:beta_i} for $i=1, 2, \cdots, N-1$ exact.  

\begin{remark}
\label{Rem:power_boundary}
If a small variation is introduced to variables $\alpha_i$ $v_i$ and $\beta_i$, it cannot be guaranteed that the power and force boundaries are not violated since $F_i$ is varied as well. The strategy is then to introduce a series of small variations and to maintain $F_n\Delta d,~n=2, 3,...i$ to remain constant except for $F_1\Delta d$. We can demonstrate that all these small variations, made possible by the strict equality of relaxation constraints \eqref{for:alpha_i} and \eqref{for:beta_i}, can eventually reduce the value of $F_1\Delta d$ and thus bring down the objective function without violating any constraints. Using the proof by contradiction, this will prove the assumed optimal solution on the strict inequality cannot be the optimal solution since small variations are possible to be obtained to further bring down the value of the objective function. 
\end{remark}

We make the notations as follows:
\begin{itemize}
    \item $\varepsilon_{i,\alpha}$ to represent a small variation for $\alpha_i$ and $\varepsilon_{i,\alpha}>0$;
    \item $\varepsilon_{i,\beta}$ to represent a small variation of the original value of $\beta_i$ and $\varepsilon_{i,\beta}>0$;
    \item $\varepsilon_{i,v}$ to represent a small variation of the original value of $v_i$ and $\varepsilon_{i,v}>0$. 
\end{itemize}

These small variables are all positive and we aim to demonstrate that there is always a set of $\varepsilon_{i,\alpha}$, $\varepsilon_{i,\beta}$ and $\varepsilon_{i,v}$ that can bring down the objective function value as long as  the solution is lying on the strict inequality of \eqref{for:alpha_i} and \eqref{for:beta_i}. We assume that all the changed variables ($v_i, \alpha_i, \beta_i$ for $i=1, 2, \cdots, N$) remain positive after the variations are applied.

\subsection{Proof of exactness: Part 1}
\label{subsec:PoE_part1}
Assume that the optimal solution $\alpha_i$ and $v_i$ for $i=1,2,\cdots,N$ exists with $\alpha_i v_i>1$. There are small variations $\varepsilon_{i,\alpha}$ and $\varepsilon_{i,v}$ that make the following inequality valid. 
\begin{equation}
    (v_i-\varepsilon_{i,v})(\alpha_i-\varepsilon_{i,\alpha})\geq 1
\end{equation}
 
\begin{remark}
In this case, $\beta_i$ remains unchanged while $v_i$ is reduced with $\varepsilon_{i,v}$. This does not violate the constraints $v_i^2 \leq \beta_i$
\end{remark}

To make sure $F_i\Delta d$ remains unchanged, we reduce the value of $\beta_{i-1}$ so that $F_i\Delta d$ remains constant. By observing  \eqref{for:ene_conserv_1_cvx}, to maintain the value of $F_i\Delta d$ unchanged, \eqref{for:epsilon_i-1_v} should be applied:

\begin{equation}
\label{for:epsilon_i-1_v}
    \varepsilon_{i-1,\beta}=\frac{\Delta dB}{0.5M}\varepsilon_{i,v} 
\end{equation}
where $\varepsilon_{i-1,\beta}$ is also a small variation related to $\varepsilon_{i,v}$.

To ensure that \eqref{for:beta_i} is valid for the $(i-1)^{th}$ segment, we introduce an equal reduction on both side of the inequality. \eqref{for:beta_i-1} is applied:
\begin{align}
\label{for:beta_i-1}
    \beta_{i-1}-\varepsilon_{i-1,\beta}&\geq v_{i-1}^2-\varepsilon_{i-1,\beta}\nonumber\\ &=(v_{i-1}-\varepsilon_{i-1,v})^2\nonumber\\
    \varepsilon_{i-1,v}&=v_{i-1}-\sqrt{v_{i-1}^2-\varepsilon_{i-1,\beta}}\nonumber\\
    \varepsilon_{i-1,\beta}&=2\varepsilon_{i-1,v}v_{i-1}-\varepsilon_{i-1,v}^2.
\end{align}

An introduction of $\varepsilon_{i,v}$ can be used to calculate $\varepsilon_{i-1,\beta}$  based on \eqref{for:epsilon_i-1_v}, and $\varepsilon_{i-1,\beta}$ can then be used to calculate a small variation on $v_{i-1}$, i.e. $\varepsilon_{i-1,v}$,  based on \eqref{for:beta_i-1}.

$\varepsilon_{i-1,v}$ and $\varepsilon_{i-1,\beta}$ are both small reduction introduced to $v_{i-1}$ and $\beta_{i-1}$ so that the model constraints are not violated. To ensure $F_{i-1}\Delta d$ to remain constant, a small reduction on $\beta_{i-2}$ needs to be introduced and it is thus defined as follows. 
\begin{equation}
    0.5M\varepsilon_{i-2,\beta}=(0.5M+C\Delta d)\varepsilon_{i-1,\beta}+B\Delta d \varepsilon_{i-1,v}
\end{equation}

This can help calculate the value of $\varepsilon_{i-2,v}$ based on a similar calculation presented in \eqref{for:beta_i-1} and this process can be iteratively continued via \eqref{for:epsilon_beta_iter} and \eqref{for:epsilon_v_iter} until the first segment is reached. 

\begin{align}
\label{for:epsilon_v_iter}
\varepsilon_{n,v}&=v_{n}-\sqrt{v_{n}^2-\varepsilon_{n,\beta}}\nonumber\\
& n=i-1, i-2, \cdots, 1
\end{align}

\begin{align}
\label{for:epsilon_beta_iter}
    0.5M\varepsilon_{n-1,\beta}&=(0.5M+C\Delta d)\varepsilon_{n,\beta}+B\Delta d \varepsilon_{n,v}\nonumber\\
    & n=i-1, i-2, \cdots, 2
\end{align}

\begin{remark}
\label{rem:small_eps}
 A sufficiently small $\varepsilon_{i,v}$ leads to a set of sufficiently small $\varepsilon_{n,v},~n=i, i-1, ..., 1$ and $\varepsilon_{n,\beta}, ~~n= i-1, ..., 1$.
 \end{remark}

Next we will demonstrate that the change of $v_n,~n=i-1, ..., 1$ will not violate the constraint defined by the inequality $1\leq v_n\alpha_n$. 

\begin{remark}
\label{rem:small_alpha}
Every decrease of $v_n, ~n=i-1, i-2, ..., 1$ leads an increase of $\alpha_n,~n=i-1, i-2, ..., 1$ to ensure $v_n\alpha_n\geq 1$ to remain valid. Due to the constraints defined by $\sum_{i=1}^N\alpha_i=T/\Delta d$, we make a special case for $v_i$ and $\alpha_i$ so that both of them are reduced. The strict inequality $1<v_i\alpha_i$ makes this possible without violating the model constraints. The reduction of $\alpha_i$ will then be used to offset the increase of $\alpha_n,~n=i-1, i-2, ..., 1$. 
\end{remark}

To calculate $\alpha_n,~ n=i-1, i-2, \cdots, 1$, we assume the relationship applies to the initial optimal solution. This is certainly valid as the initial optimal solution should be a feasible solution for Model \eqref{for:EETC_cvx}. 

\begin{equation}
    v_n \alpha_n=D_n\geq 1
\end{equation}

After the introduction of a small variation, the product remains constant defined as:
\begin{equation*}
    (v_n-\varepsilon_{n,v})(\alpha_n+\varepsilon_{n,\alpha})=D_n
\end{equation*}
where, $n= i-1, i-2,\cdots,1$. 

Thus, we can obtain the value of $\varepsilon_{n,\alpha},~n= i-1, i-2,\cdots,1$ by
\begin{equation}
\label{for:epsilon_n_alpha}
    \varepsilon_{n,\alpha}=\frac{\alpha_n\varepsilon_{n,v}}{v_n-\varepsilon_{n,v}}. 
\end{equation}

\begin{remark}
\label{rem:mono}
Based on Remarks \ref{rem:small_eps} and \ref{rem:small_alpha}, we can use the first derivative to evaluate the relationship between a small variable $v_n$ and final total increase of $\alpha_n,~ n=1, 2, \cdots, i-1$. 
\end{remark}

According to \eqref{for:beta_i-1}, the derivative of $\varepsilon_{i-1,v}$ over $\varepsilon_{i-1,\beta}$ is given as follows:
\[
\frac{d (\varepsilon_{i-1,v})}{d(\varepsilon_{i-1,\beta})}=\frac{1}{\sqrt{v_{i-1}^2-\varepsilon_{i-1,\beta}}}>0
\]

By observing \eqref{for:epsilon_beta_iter}, it is known that this is a strictly increasing monotonic function over either $\varepsilon_{i-1,\beta}$ and $\varepsilon_{i-1,v}$

By observing \eqref{for:epsilon_n_alpha}, we are able to obtain its first derivative of $\varepsilon_{n,\alpha}$ over $\varepsilon_{n,v}$ as follows:
\[
\frac{d (\varepsilon_{n,\alpha})}{d(\varepsilon_{n,v})}=\frac{\alpha_n}{v_n-\varepsilon_{n,v}}+\frac{\alpha_n\varepsilon_{n,v}}{(v_n-\varepsilon_{n,v})^2}=\frac{D_n}{(v_n-\varepsilon_{n,v})}>0
\]

The reduction of $\alpha_i$ will be equal to the sum of all $\varepsilon_{n,\alpha},~n=i-1, i-2, \cdots, 1$, represented by:
\begin{equation}
\label{for:epsilon_i_alpha}
    \varepsilon_{i,\alpha}=\sum_{n=1}^{i-1}\varepsilon_{n,\alpha}.
\end{equation}

Given that each $\varepsilon_{n,\alpha}$ is a strictly increasing monotonic function over $\varepsilon_{n,v}$ which is also a strictly increasing monotonic function over $\varepsilon_{i,v}$, by the property of monotonic functions, it can be deducted that $\varepsilon_{n,\alpha}$ is a strictly increasing monotonic function over $\varepsilon_{i,v}$. 

Therefore, there exist sufficiently small numbers $\varepsilon_{n,v}$ and $\varepsilon_{n,\alpha}$ that reduce the value of the assumed optimal solution $v_i$ and $\alpha_i$ which are constrained by a strict inequality $1<\alpha_iv_i$ to generate a new feasible solution $\tilde{v}_n$ and $\tilde{\alpha}_n$, defined by \eqref{for:tilde_v_n} and \eqref{for:tilde_alpha_n} respectively: 

\begin{equation}
\label{for:tilde_v_n}
    \tilde{v}_n=v_n-\varepsilon_{n,v},~ n= i, i-2,\cdots,1
\end{equation}

\begin{equation}
\label{for:tilde_alpha_n}
    \tilde{\alpha}_n=
    \begin{cases}
    \alpha_n-\varepsilon_{n,\alpha} \text{ if }~n= i\\
    \alpha_n+\varepsilon_{n,\alpha} \text{ if }~n= i-1, i-2, \cdots, 1
    \end{cases}
\end{equation}

\begin{figure}[!th]
    \centering
    \includegraphics[width=0.8\textwidth]{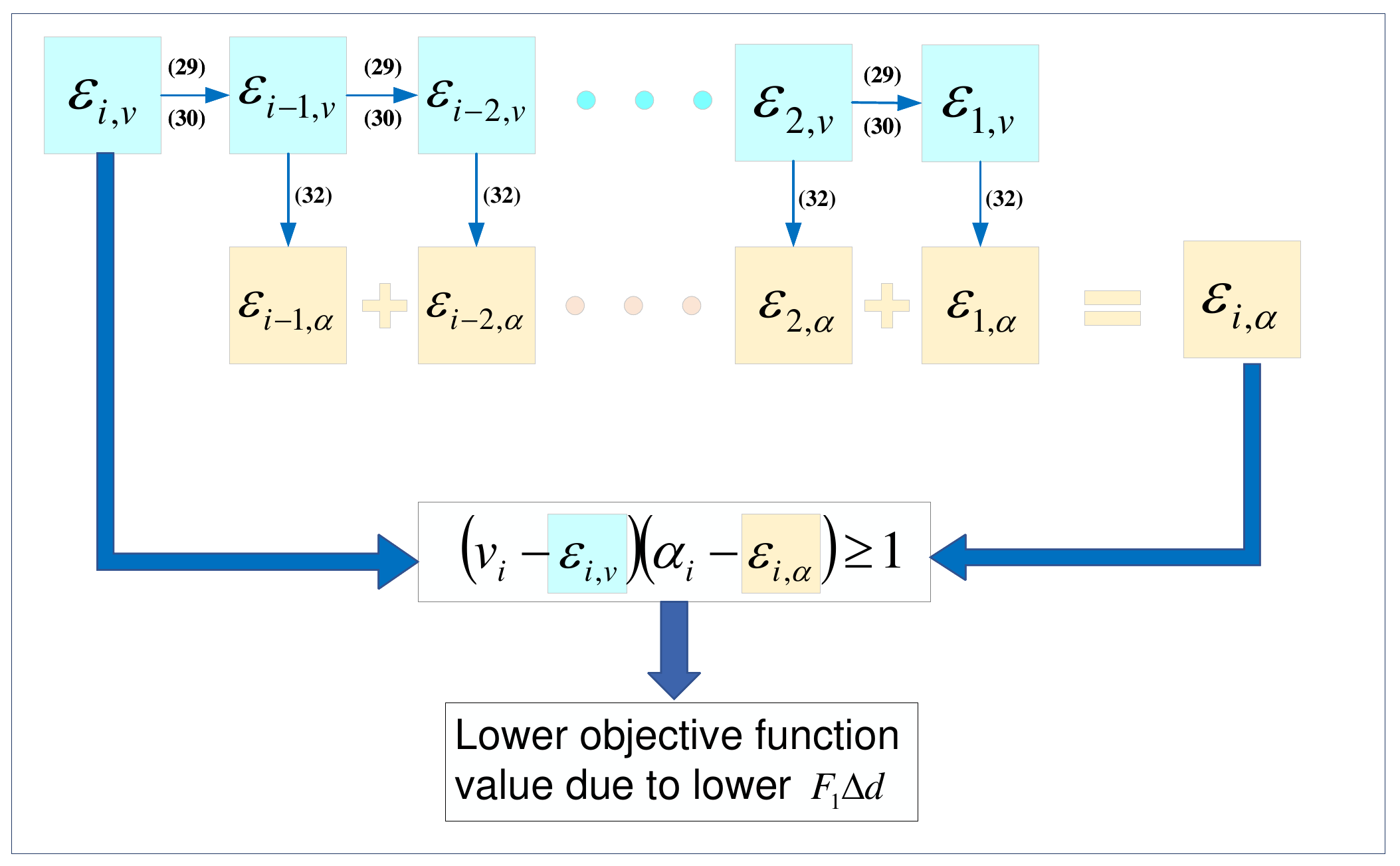}
    \caption{An illustrative graph for the part 1 of exactness proof: a small variation of $v_i$ leads to a small variation of $\alpha_i$ without violating constraints and finally reduce the objective function value. Relevant applied equation number for calculations is listed in the figure.}
    \label{fig:PoE_Part_1}
\end{figure}

Fig. \ref{fig:PoE_Part_1} shows an illustration for the exactness proof in part 1. On the segments other than the first one, $F_i\Delta d$ remains constant. On the first segment, the initial value of $\beta_1$ and $v_1$ will be reduced with a small value $\varepsilon_{1,\beta}$ and $\varepsilon_{1,v}$ and this does not violate the power and force boundaries and other constraints as claimed in Assumption \ref{assum:1}.

In other words, the new set of variables is the feasible solution and  attains a lower value of the objective function. This contradicts the initial assumption that $v_i$ and $\alpha_i$ are the optimal solution. This shows that there is no such a solution $v_i$ and $\alpha_i$ that is on the strict inequality $1<\alpha_i v_i$ and also achieve the optimal solution for Model A. Model B will always achieve its optimal solution on equality constraint $1=\alpha_i v_i$. This completes Part 1 of the proof of exactness. 

\subsection{Proof of exactness: Part 2}
\label{subsec:PoE_part2}
Similarly, we seek similar sets of small variations which always prevent the variables on the strict inequality from being optimal. First, assume that the optimal solution exists on the strict inequality of \eqref{for:beta_i}, and thus there are two small variations that make following constraint valid: 
\[
(v_i+\varepsilon_{i,v})^2\leq \beta_i-\varepsilon_{i,\beta}
\]

\begin{figure}[!th]
    \centering
    \includegraphics[width=\textwidth]{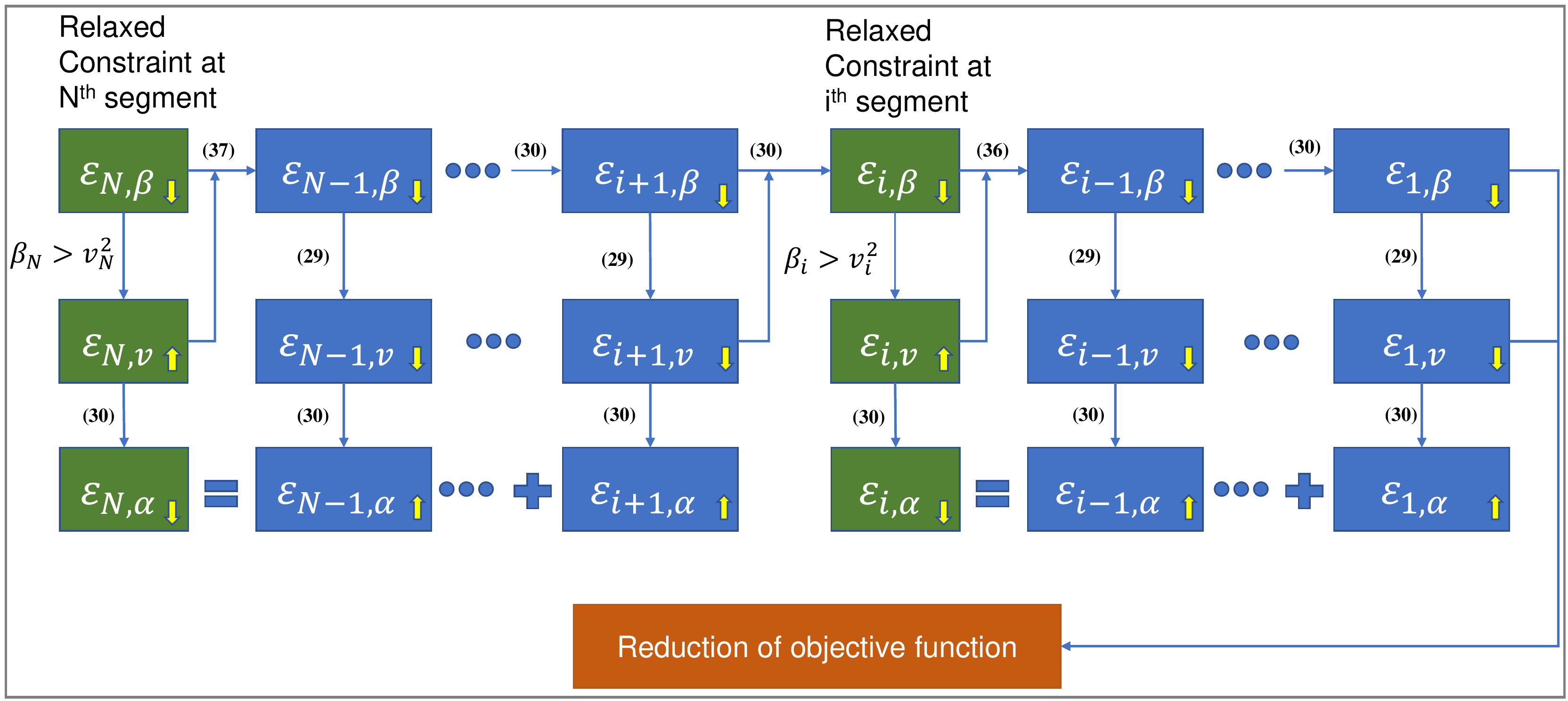}
    \caption{An illustrative graph for the part 2 of exactness proof: with an assumed relaxed constraints at the $N^{th}$ segment, we are able to achieve a feasible solution using a set of small variations and finally reduce the value of the objective function. The yellow arrow within each block shows if the variation is to increase or to decrease the original variable and the relevant applied equation for calculation of different variables is included from one variable block to another.}
    \label{fig:PoE_Part_2}
\end{figure}

In Fig. \ref{fig:PoE_Part_2}, we illustrate the part 2 of exactness proof for constraint \ref{for:beta_i}. 

We introduce a small variation $\varepsilon_{i-1,\beta}$ to bring down the initial optimal speed $v_{i-1}$ and it can be demonstrated that using \eqref{for:epsilon_v_iter} and \eqref{for:epsilon_beta_iter}, we can bring down the value of initial value of $\beta_n$ and $v_n$ where $n=i-1, i-2, \cdots,1$ without violating the constraints. 

In addition, by adopting \eqref{for:epsilon_n_alpha}, it is possible to calculate the increase of $\alpha_n,~ n=i-1, i-2, \cdots, 1$ and the value of the reduction of $\alpha_i$ can be calculated by \eqref{for:epsilon_i_alpha}. Based on \eqref{for:epsilon_i_alpha}, it is possible to calculate the small increase of $v_i$, i.e. $\varepsilon_{i,v}$. To maintain a constant value of $F_i\Delta d$, based on \eqref{for:ene_conserv_1_cvx}, the following two relations need to be applied. 

\begin{equation}
\label{for:eps_beta_i}
    \varepsilon_{i,\beta}(C\Delta d+0.5M)=B\Delta d\varepsilon_{i,v}+0.5M\varepsilon_{i-1,\beta}
\end{equation}

\begin{equation}
\label{for:eps_beta_N}
    \varepsilon_{N,\beta}(C\Delta d+0.5M)=B\Delta d\varepsilon_{N,v}+0.5M\varepsilon_{N-1,\beta}
\end{equation}

Based on Assumption \ref{assum:2}, there exists a small gap between $\beta_N$ and $v_N$ and a small variation can be applied to reduce the value of $\beta_N$ and $\beta_N-1$ and increase the value of $v_N$. 

According to \eqref{for:epsilon_v_iter} and \eqref{for:epsilon_beta_iter}, a small variation $\varepsilon_{N-1,\beta}$ on $\beta_{N-1}$ gives rise to a series of reduction on $v_n,~n=N-1, N-2, \cdots, i+1$ and based on Remark \ref{rem:mono}, it is known that $\sum_{n=i+1}^{N-1}\alpha_n$ is strictly increasing monotonic over $\varepsilon_{N-1,\beta}$. Since the speed $v_n$ is reducing, $\alpha_n$ is increasing. On the other hand, $\varepsilon_{N,\alpha}$ is a strictly increasing monotonic function of $\varepsilon_{N,v}$. Note that $v_N$ is increasing with a variation of $\varepsilon_{N,v}$. 

Now define a function $f(\varepsilon_{N-1,\beta})$ as follows:
\begin{equation}
\label{for:fun_eps_beta_N-1}
    f(\varepsilon_{N-1,\beta})=\sum_{n=i+1}^{N-1}\alpha_n-\varepsilon_{N,\alpha}
\end{equation}

This function is a continuous function since no extreme value is taken and has a feasible domain of $[0, \bar{\varepsilon}_{N-1,\beta}]$ where $\bar{\varepsilon}_{N-1,\beta}$ is determined by \eqref{for:eps_beta_N} when $\varepsilon_{v,N}$ takes a value of zero, as shown below:
\[ 
 \bar{\varepsilon}_{N-1,\beta}=\left(1+\frac{\Delta d C}{0.5M}\right)\varepsilon_{N,\beta}. 
\]

On the boundary of the feasible domain, \eqref{for:fun_eps_beta_N-1} takes values of $-\varepsilon_{N,\alpha}<0$ and $\sum_{n=i+1}^{N-1}\alpha_n>0$. Based on \textbf{Intermediate Value Theorem}, there exists one value of $\varepsilon_{N-1,\beta}$ that causes the function to take a value of zero. We use $\tilde{\varepsilon}_{N-1,\beta}$ and $\tilde{\varepsilon}_{i,\beta}$ to denote the corresponding values that leads \eqref{for:fun_eps_beta_N-1} to take a value of zero. With the existence of $\tilde{\varepsilon}_{N-1,\beta}$ and $\tilde{\varepsilon}_{i,\beta}$, it is known that every $F_n\Delta d$ for $n=i+1, i+2, \cdots, N$ will be maintained constant and the constraints on $\alpha_n$  for $n=i+1, i+2, \cdots, N$ remain valid as the increase of $\alpha_n$ for $n=i+1, i+2, \cdots, N-1$  is offset by the decrease of $\alpha_{N}$. 

Given a sufficiently small variation of  $\tilde{\varepsilon}_{i,\beta}$, we are able to conduct similar analysis based on \textbf{Intermediate Value Theorem} and \eqref{for:eps_beta_i} to show that there exists $\tilde{\varepsilon}_{i-1,\beta}$ and $\tilde{\varepsilon}_{i,v}$ that makes the constraints on $\alpha_n~~n=1,2, \cdots, i$ valid and meanwhile brings down a series of reduction on $v_n$ and $\beta_n$ for $n=1,2, \cdots, i-1$. Based on Assumption \ref{assum:1}, the existence of variations will decrease the objective function on the first segment and this proves that \eqref{for:beta_i} remain exact for any index number $n=1,2, \cdots, N-1$. \textbf{It is worth noting that for the case of the $N^{th}$ segment, the $\beta_N$ will always keep approaching $v_N^2$ and virtually this constraint remains exact even with the assumption on its strict inequality in Assumption \ref{assum:2}}.

In a summary, we have demonstrated that all optimal solution of Model B will only be obtained on the equality of relaxation constraints defined by \eqref{for:alpha_i} and \eqref{for:beta_i}. This thus proves that Theorem \ref{the:exactness} is true based on Assumption \ref{assum:1} and \ref{assum:2}. As a matter of fact, as long as there is any $F_n\Delta d$ where $n=1, 2, \cdots i$ where a small reduction of $F_n\Delta d$ does not violate the power and force boundary defined by \eqref{for:power_boundary_cvx} and \eqref{for:force_boundary}, e.g. $F_n\Delta d>0$, the proof remains valid.

\section{Numerical verification on a case with a steep gradient}
\label{sec:num_veri}
In this section, we demonstrate a quick numerical verification on the model exactness using a case study with a steep gradient. 

The modeling parameters are listed in Table \ref{tab:para_model}.
\begin{table}[th]
    \caption{Modeling parameters of the traction system}
    \label{tab:para_model}
    \centering
    \begin{tabular}{p{8cm} p{2.5cm}<{\centering}}
    \toprule
    Parameter & Value \\
    \midrule
    $T(\rm s)$ & 260\\
    $M (\rm t)$ & 144   \\
    $F_{t,max} (\rm kN)$ & 230.81    \\
    $F_{b,max} (\rm kN)$ & 230.81    \\
    $P_{t,max} (\rm kW)$ & 2520  \\
    $P_{b,max} (\rm kW)$ & 2520  \\
    $\eta_t$ & 0.9  \\
    $\eta_b$ & 0.6  \\
    $A ( \rm kN)$           & 3.0016    \\
    $B (\rm kN/(km/h))$     & 2.016e-2  \\
    $C (\rm kN/(km^2/h^2))$ & 6.9692e-4 \\
    \bottomrule
    \end{tabular}
\end{table}

\begin{figure*}[!th]
    \centering
    \includegraphics{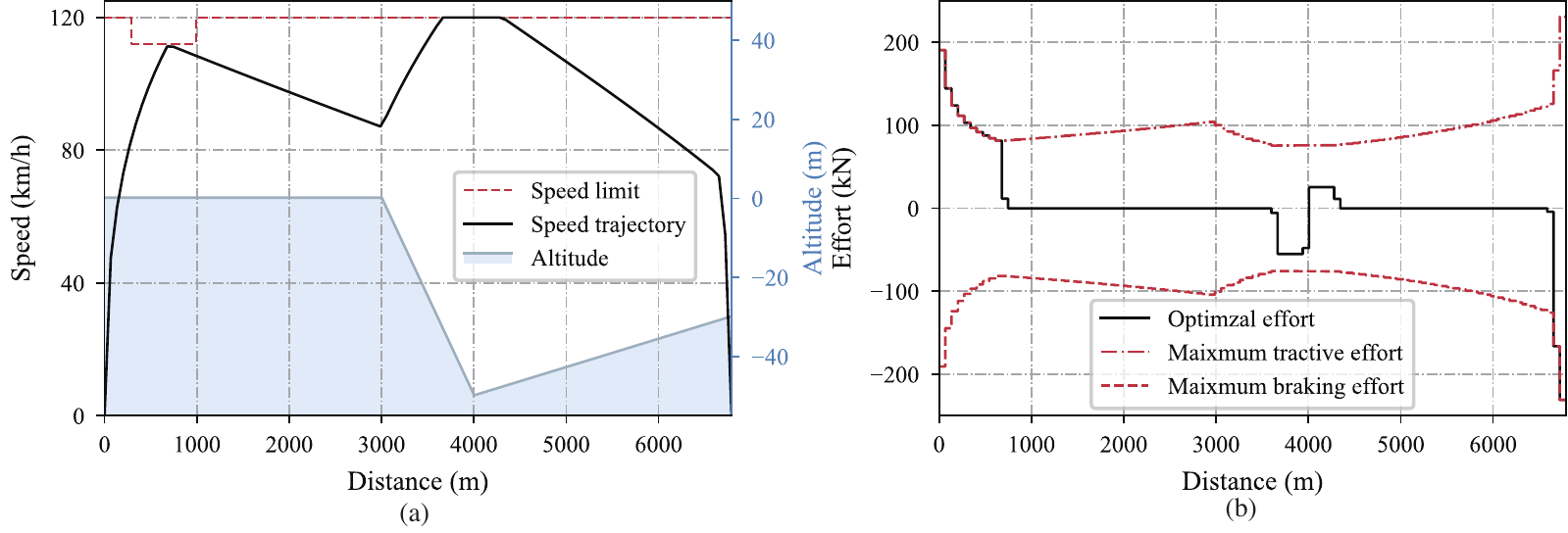}
    \caption{Optimal speed trajectory and tractive/braking efforts for a typical urban rail case based on convex optimization. (a) It demonstrates the optimal speed trajectory with varying speed limit and steep gradient. (b) The tractive/braking effort shows that the train optimal operation consists of maximum traction/braking, partial traction/braking, and coasting.}
    \label{fig:speed_effort}
\end{figure*}

In Fig. \ref{fig:speed_effort}, we present a typical urban rail case with an arbitrary steep downhill gradient between 3000 m and 4000 m on which the train needs to brake to maintain a constant speed. On distance 4000 m, the change of the gradient has been counteracted by the varying efforts. The results show that the model can be solved and quickly achieve an optimal solution with a steep gradient.

\begin{figure*}[!th]
    \centering
    \includegraphics{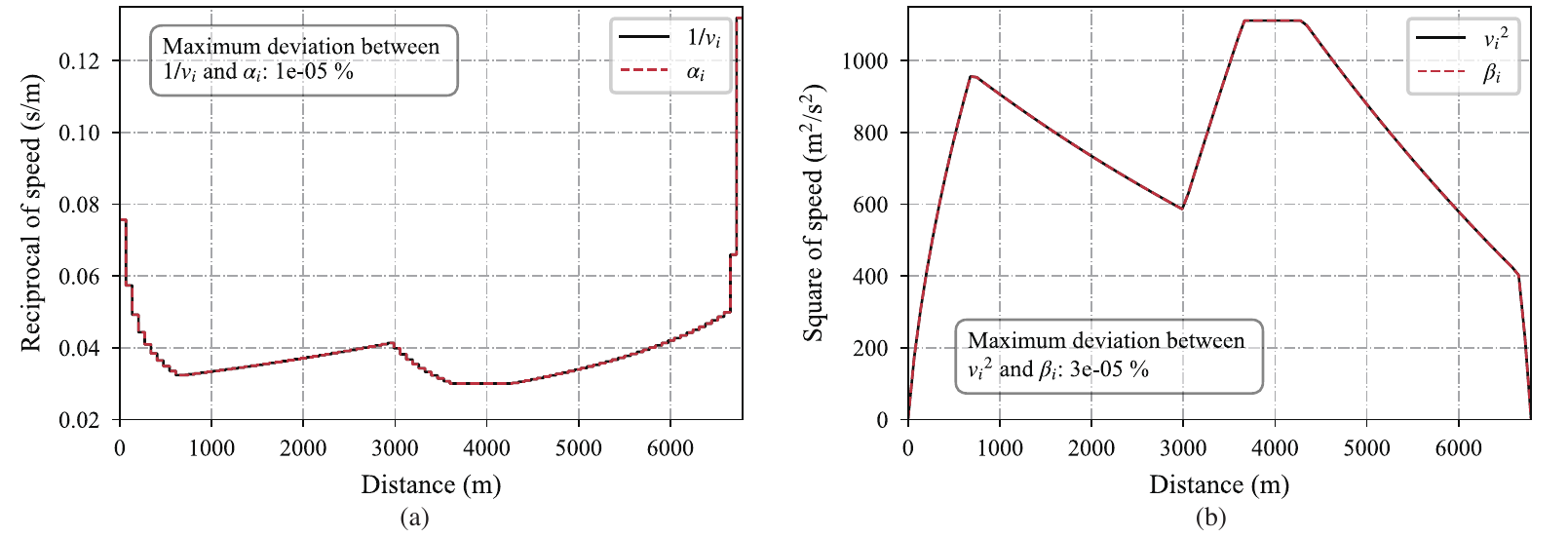}
    \caption{Comparison of variables. (a) Comparison between the reciprocal of the candidate speed $1/v_i$ and  $\alpha_i$. (b) Comparison between the square of the candidate speed $v_i^2$ and the variable $\beta_i$. Both figures show that when the optimization model achieves the optimal solution, inequality constraints \eqref{for:alpha_i} and \eqref{for:beta_i} achieves equality.}
    \label{fig:aux_var}
\end{figure*}

Fig. \ref{fig:aux_var} shows a comparison of two sets of variables of interest and demonstrates that equality for \eqref{for:alpha_i} and \eqref{for:beta_i} will be attained when the optimal solution is achieved. More numerical results and discussion can be found in our newly published preprint paper\cite{FengM2022convex}.

\section{Conclusions}
\label{sec:conclusion}
As a continued work of \cite{FengM2022convex}, this paper presents the exactness proof for an EETC model based on convex optimization referred to as Model B. We applied two constraint relaxation to convert an non-convex Model, i.e. Model A, into a convex Model B. The convexity of Model B was first verified and the exactness proof of Model B was then conducted. With its convexity and relaxation exactness, the proposed model with mature solution approaches such as interior-point barrier methods boast a very high computational efficiency and global optimality compared to its counterpart such as MILP, Psuedospectral Method and other non-linear programming methods. Its high flexibility could enhance the model's adaptivity for other types of transportation modes to further reduce the energy consumption.

%\clearpage
\bibliographystyle{unsrt}
\bibliography{Resources/bibliography}
\end{document}